\documentclass[a4paper]{article}
\usepackage{amssymb,graphicx,url,amsmath,subfigure
}
\usepackage{times}
\usepackage{tikz}
\usepackage{soul}
\usepackage{color}
\usepackage{xcolor}
\usepackage{algorithm}
\usepackage{algorithmic}
\usepackage{enumerate}
\usepackage{a4wide}
\usepackage{comment}
\usepackage{authblk}

\newcommand{\Exp}{\mathbf{E}}

\newcommand{\R}{{\mathbb R}}
\newcommand{\N}{{\mathbb N}}
\newcommand{\C}{{\mathbb C}}

\newcommand{\cM}{{\cal M}}
\newcommand{\cB}{{\cal B}}

\newtheorem{Theorem}{Theorem}
\newtheorem{Lemma}{Lemma}
\newtheorem{Definition}{Definition}

\date{July 21, 2017}

\title{A central limit like theorem
for Fourier sums}
\author[1]{Dominik Janzing
\footnote{corresponding author: dominik.janzing@tuebingen.mpg.de}}
\author[2]{Naji Shajarisales}
\author[1,3]{Michel Besserve}

\affil{Max Planck Institute for Intelligent Systems, T\"ubingen, Germany}
\affil[2]{Carnegie Mellon University, Pittsburgh, USA}
\affil[3]{Max Planck Institute for Biological Cybernetics, T\"ubingen, Germany}

\begin{document}
\maketitle

\begin{abstract}
We consider the probability distributions of
values in the complex plane attained by Fourier sums of the form
$\frac{1}{\sqrt{n}}\sum_{j=1}^n 
a_j e^{-2\pi i j \nu}$ when the frequency $\nu$ is drawn uniformly at random from 
an interval of length $1$. 
If the coefficients $a_j$ are i.i.d.~drawn with finite third moment, 
the distance of these distributions to an isotropic two-dimensional
Gaussian on $\C$ converges in probability to zero for any pseudometric
on the set of distributions for which the distance between empirical distributions
and the underlying distribution converges to zero in probability.
\end{abstract}

The classical version of the central limit theorem states that for 
a series of real-valued  independent identically distributed (iid) random variables $X_1,X_2,\dots$ with 
$\Exp[X_j]=0$ and finite variance $\sigma^2$ the sequence 
\begin{equation}
\frac{1}{\sqrt{n}} \sum_{j=1}^n X_j
\end{equation}
converges in distribution to a Gaussian random variable with zero mean
and variance $\sigma^2$ \cite[p.357]{billingsley1995}. 
Formulations for random vectors $X_j$ state convergence to multi-variate Gaussians
\cite{vaart1998}. 
Other well-known generalizations drop the assumption `identically distributed'
and replace it, for instance, with the Lyapunov condition 
\[
\lim_{n\to \infty} \frac{1}{s_n^{2+\delta}} \sum_{i=1}^n \Exp[X_j^{2+\delta}] =0,
\]
for some $\delta>0$, where $s_n$ denotes the sum of all variances of $X_1,\dots,X_n$
\cite[p.362]{billingsley1995}.
Then 
\[
\frac{1}{\sqrt{s_n}} \sum_{j=1}^n X_j
\]
converges in distribution to a standard Gaussian. 
It is also known that the independence assumption can be replaced with
appropriate notions of weak dependence, e.g. \cite[Theorem~27.4]{billingsley1995}.
However, significantly more general scenarios yield Gaussians as limiting distributions. Here we consider sequences of Fourier sums of the form
\begin{equation}\label{eq:ahat}
\hat{a}^n(\nu) :=\frac{1}{\sqrt{n}} \sum_{j=1}^n 
a_j   e^{-2\pi i \nu j},
\end{equation}
and show that 
sampling from random frequencies yields asymptotically a Gaussian -- in a sense to be specified below -- if the coefficients $a_j$ are i.i.d.~drawn.
More precisely, let $A_1,A_2,\dots$ be a sequence of real-valued i.i.d.~variables 
on a probability space $(\Omega,\Sigma,P_\Omega)$. Then we
first define for each frequency $\nu$ the sequence $(\hat{A}^n_\nu)_{n\in \N}$ of random variables 
via 
\[
\hat{A}^n_\nu := \frac{1}{\sqrt{n}} \sum_{j=1}^n 
A_j   e^{-2\pi i \nu j}.
\]
Using known vector-valued central limit theorems
one can easily show that, under some technical condition of Lyapunov type  detailed below, $\hat{A}^n_\nu$ converges
to a Gaussian on the complex plane for each $\nu$ since it is obtained by a sum of the independent
(but not identically distributed)
complex-valued random variables 
$X_j:=A_j e^{-2\pi i \nu j}$.

Here, however, we define for each $\omega \in \Omega$
the sequence $(\hat{A}^n_.(\omega))_{n\in \N}$ of random variables 
on the probability space $([-1/2,1/2],\cB,\lambda)$, with $\cB$ denoting the Borel sigma algebra and
$\lambda$ the Lebesque measure (to formalize the random choice of a frequency), via
\[
\hat{A}^n_. (\omega) :\quad  \nu \mapsto  \hat{A}^n_\nu(\omega)= \emph{•}\frac{1}{\sqrt{n}}
\sum_{j=1}^n A_j(\omega) e^{-2\pi i \nu j}. 
\]
 
The problem is motivated by Ref.~\cite{shajarisales2015telling} which considers 
linear time invariant filters whose coefficients are randomly chosen. The question arising there was how the filter's frequency response \eqref{eq:ahat}
behaves in the limit $n\to\infty$ for  `typical' choices of filter coefficients $a_j$ when the latter are randomly drawn.
 
Note that the problem would become simple if we were to consider 
$\hat{a}^n(\nu)$ at the discrete frequencies $\nu_l:=j/n$ for $l=1,\dots,n/2$
(in signal processing $|\hat{a}^n(\nu_l)|^2$ is also known as the periodogram of a signal  \cite{fay2001})
and $a_j$ were assumed to be drawn from independent Gaussians. For Gaussian $A_j$, the random variables 
$\hat{A}^n_{\nu_l}$
defined on the probability space $(\Omega,\Sigma,P_\Omega)$
 are also independent Gaussians for these different discrete frequencies, which can easily checked 
 by computing the covariances. 
The question changes drastically when we consider 
the full continuum of frequencies, because $\hat{A}^n_\nu$ and $\hat{A}^n_{\nu'}$ are not in general independent for $\nu \neq \nu'$. We will show, however, that they become asymptotically independent, which then results
in an appropriate limit theorem
for $\hat{A}^n_.(\omega)$.

Crucial to phrase our limit theorem is the following type of distance measures on probability distributions:
\begin{Definition}[well-behaved pseudometric]\label{def:well} 
For an arbitrary sequence $Z_1,Z_2,\dots$ of i.i.d.~random variables
on the probability space $(\Omega',\Sigma',P_{\Omega'})$ let
$P_Z$ denote the distribution of each $Z_j$ and 
$\hat{P}_{Z_1(\omega'),\dots,Z_k(\omega')}$
denote the empirical distribution after the first $k$ samples. 

Let ${\cal M}_l$ denote the set of probability measures on the Borel-measurable subsets of $\R^l$.
Then a pseudometric $d:\cM_l \times \cM_l\rightarrow \R^+_0$ is called `well-behaved' if the distance between $P_Z$ and 
$\hat{P}_{Z_1(\omega'),\dots,Z_k(\omega')}$ 
converges in probability to zero uniformly over all i.i.d.~sequences. More precisely, 
for every $\epsilon,\delta>0$ there is a $k_0$ such that for all $k\geq k_0$
\[
P_{\Omega'} \left\{ d(P_Z,\hat{P}_{Z_1(\omega'),\dots,Z_k(\omega')}) \geq \epsilon \right\} \leq \delta,
\]
holds for all sequences $Z_1,Z_2,\dots$ and probability spaces $(\Omega',\Sigma',P_{\Omega'})$. 
\end{Definition}
For distributions $Q,R$ on $\R$ 
with cumulative distribution functions
$F_Q$ and $F_R$, respectively, $d(Q,R):=\|F_Q-F_R\|_\infty$ provides
a simple example of a well-behaved distance since
\[
\|F_{\hat{Q}_k} -F_Q\|_\infty \geq \epsilon,
\]
occurs with probability at most $2e^{-2k \epsilon^2}$ due to
Massart's formulation \cite{Massart90} of the  Dvoretzky-Kiefer-Wolfowith (DKW) inequality.
Another example is given by $d(Q,R):=\sup_B |Q(B)-R(B)|$ where $B$
runs over some set of sets whose indicator functions have finite VC-dimension.
This follows from Vapnik and Chervonenkis' uniform bound on the deviation
of empirical frequencies of events from the corresponding probabilities \cite{vapnik98}.
Using
Reproducing Kernel Hilbert Spaces (RKHS) one can construct a further example:
the so-called kernel mean embedding \cite{smola2007} represents distributions as vectors
in a Hilbert space. Then the Hilbert space distance is a well-behaved metric.
This follows easily from the uniform consistency result in \cite[Theorem~4]{gretton2006}
for the empirical estimator of this distance\footnote{If both samples in 
\cite[Theorem~4]{gretton2006}
 are drawn from the same distribution and one of the sample sizes tends to infinity the 
 bound describes the distance between empirical and true distribution.}. 

The purpose of this article is to show the following result:
\begin{Theorem}\label{thm:main}
Let $P_{\hat{A}^n_.(\omega)}$ denote the distribution  of $\hat{A}^n_.(\omega)$ and
$G$ the distribution on $\C$ for which
 real and imaginary parts are independent
 Gaussians with
mean zero and variance $1/2$. Then the distance between
 $P_{\hat{A}^n_.(\omega)}$ and $G$ converges
 to zero in probability for every well-behaved pseudometric $d$. More, precisely, the random variable 
\[ 
\omega \mapsto d(P_{\hat{A}^n_.(\omega)},G)
\] 
converges to zero in probability.
\end{Theorem}
Obviously,
the interval $[-1/2,1/2]$ can be replaced by any interval of length $1$, as stated in the abstract. 

The first step of the proof will be to investigate the asymptotics of
the variances and covariances of real and imaginary part of $\hat{A}^n_\nu$ and the covariances between real and imaginary parts
$\hat{A}^n_\nu$ and $\hat{A}^n_{\nu'}$ for different frequencies $\nu,\nu'$. To this end,
we represent complex numbers as vectors in $\R^2$ and obtain the following result:
\begin{Lemma}[asymptotic covariances]\label{lem:covasym}
Let $\nu_1,\dots,\nu_k$ be some  arbitrary non-zero frequencies in $(-1/2,1/2)$ 
with $|\nu_i|\neq |\nu_j|$ for $i\neq j$. 
Let $C^{(n)}$ denote the covariance matrix of the 
random vector $\frac{1}{\sqrt{n}}\sum_{j=1}^n S_j$ with
\[
S_j := A_j 
\Big(\cos(2\pi \nu_1 j),\sin(2\pi \nu_1 j),
\cos(2\pi \nu_2 j),\sin(2\pi \nu_2 j),\dots,
\cos(2\pi \nu_k j),\sin(2\pi \nu_k j)\Big)^T.
\]
Then 
\[
\lim_{n\to\infty} C^{(n)} = \frac{1}{2} {\bf 1},
\]
where ${\bf 1}$ denotes the identity in $2k$ dimensions.
\end{Lemma}

\noindent
Proof: We first introduce the vector 
$c:=(1,0)^T$ and
the rotation matrix
\[
D_\nu := \left(\begin{array}{cc} \cos(2\pi \nu ) & -\sin (2 \pi \nu) \\ \sin (2\pi \nu ) & \cos (2\pi \nu) \end{array} \right).
\]
Using powers of these rotations, we can write 
the random vector $S_j$ as the direct sum
\[
S_j := A_j  \left[ D^j_{\nu_1} c \oplus  D^j_{\nu_2} c \oplus \cdots \oplus D^j_{\nu_k} c\right].  
\]
Its covariance matrix reads 
\[
C_j := \underbrace{\left(\begin{array}{cccc}
D_{\nu_1}^j & & & \\
& D_{\nu_2}^j & & \\
& & \ddots & \\
& & & D_{\nu_k}^j \end{array}\right)}_{D^j}
\underbrace{ 
\left(\begin{array}{cccc}
c c^T & cc^T & \cdots & ccT \\
cc^T  &  \ddots &    & \vdots \\    
\vdots &         &  &       \\ 
cc^T & \cdots & & cc^T 
\end{array}\right)
}_{C_0}
\underbrace{
\left(\begin{array}{cccc}
D_{\nu_1}^{-j} & & & \\
& D_{\nu_2}^{-j} & & \\
& & \ddots & \\
& & & D_{\nu_k}^{-j} \end{array}\right)}_{D^{-j}}.  
\]
Since the random vectors $S_1,\dots,S_n$ are uncorrelated (because the variables $A_j$ are independent and thus uncorrelated), the weighted sum 
\[
S^{(n)}:=\frac{1}{\sqrt{n}} \sum_{j=1}^n S_j
\]
has the covariance matrix 
\[
C^{(n)} 
:=\frac{1}{n} \sum_{j=1}^n D^j C_0 D^{-j}. 
\]
Block $ll'$ within the $k\times k$ block matrices of format $2\times 2$
reads
\[
C^{(n)}_{ll'} 
= \frac{1}{n} \sum_{j=1}^n F_{ll'}^j (cc^T),
\]
where $F^j_{ll'}$ denotes the $j$th power of the map
$F_{ll'}$ on the space ${\cal M}_2(\C)$ of complex-valued $2\times 2$-matrices defined via 
\[
F_{ll'}(M):= D_{\nu_l} M D^{-1}_{\nu_{l'}}. 
\]
Note that $F_{ll'}$ is a unitary map on ${\cal M}_2(\C)$ with respect to the inner product $\langle A, B\rangle:={\rm tr}(B^\dagger A)$, where $\dagger$ denotes the Hermitian conjugate.
Therefore, von Neumann's mean ergodic theorem \cite{ReedSimon80} implies
\[
\frac{1}{n} \sum_{j=1}^n F_{ll'}^j (cc^T) = Q_{ll'} (cc^T),
\]
where $Q_{ll'}$ denotes the orthogonal projection\footnote{We could have also applied the mean ergodic theorem to unitary map $C_0\mapsto 
D C_0 D^{-1}$ instead of applying it to each block separately, but finally this would not have simplified the analysis.}
onto the $F_{ll'}$-invariant subspace of ${\cal M}_2(\C)$.  If $r_1:=(1,i)^T$ and $r_2:=(1,-i)^T$
denote the joint eigenvectors of all $D_\nu$ with eigenvalues
 $e^{\pm 2\pi \nu}$ 
 then $F_{ll'}$ has the $4$ eigenvectors $r_{j} r_l^\dagger $
 with $j,l=1,2$ and eigenvalues 
$e^{i 2\pi (\pm \nu_l \pm \nu_{l'})},
e^{i 2\pi (\pm \nu_l \mp \nu_{l'})}
$. For $l\neq l'$ the $F_{ll'}$-invariant subspace is $0$ because
all eigenvalues differ from $1$ due to $0\neq |\nu_l|\neq |\nu_{l'}|\neq 0$. Hence, the non-diagonal blocks 
of $C^{(n)}$ vanish in the limit. 
To consider the diagonal blocks, note that
$F_{ll}$ is then  just the adjoint map of $D_{\nu_l}$, which can be restricted to
the space of real-valued symmetric matrices ${\cal M}^{sym}_2(\R)$. 
We then conclude
\[
\lim_{n \to \infty} \frac{1}{n} \sum_{j=1}^n  F^j_{ll} (M) = Q_{ll}(M) = \frac{1}{2} {\rm tr} (M) {\bf 1}
\quad \forall M \in {\cal M}^{sym}_2(\R).
\]
This is since
multiples of the identity are the only real symmetric matrices that commute
with $D_\nu$ for $\nu \in (-1/2,1/2) \setminus\{0\}$ and because $F_{ll}$ preserves the trace.
$\Box$

\vspace{0.3cm}
\noindent
We now state
the following central limit theorem for random vectors 
\cite{bentkus05}:
\begin{Lemma}[CLT for random vectors with explicit bound]\label{lem:bentkus}
Let $X_1,\dots,X_n$ be independent random vectors
in $\R^d$ such that $\Exp[X_j]=0$ for all $j$. 
Write $S:=\sum_{j=1}^n X_j$ and assume that
the covariance matrix $C_S$ of $S$ is invertible. Let $Z$ be a centered Gaussian random vector with covariance matrix $C_S$. Let ${\cal C}$ denote the set of convex subsets of
$\R^d$. Then,  
\begin{equation}\label{eq:boundbentkus}
\sup_{B \in {\cal C}} \left|P\{ S\in B \} - P\{Z \in B\}\right|  \leq \eta d^{1/4} \sum_{j=1}^n \beta_j,
\end{equation}
with
\[
\beta_j:= \Exp\left[\left\| \sqrt{C_S}^{-1} X_j\right\|^3\right],
\] 
for some $\eta>0$. 
\end{Lemma}
Since we will not use the explicit bound we derive a simpler asymptotic statement
as implication:
\begin{Lemma}[simplified CLT for random vectors]\label{lem:derived} 
Let $Y_j$ with $j\in \N^*$ independent random vectors with covariance matrices $C_j$ such that $C^{(n)}:=\frac{1}{n} \sum_{j=1}^n C_j$ 
converges to some invertible matrix $C$ with respect to any matrix norm. Assume, moreover, that there exists a constant $b<\infty$ such that $\Exp[\|Y_j\|^3] \leq b$ for all $j$.  
Then
\[
S^{(n)}:=\frac{1}{\sqrt{n}} \sum_{j=1}^n Y_j
\]
converges in distribution to a multivariate Gaussian with covariance matrix $C$. 
\end{Lemma}

\noindent
Proof: Applying Lemma~\ref{lem:bentkus}
to the variables $X_j:=\frac{1}{\sqrt{n}} Y_j$ 
yields
\[
\beta_j = 
\frac{1}{n^{3/2}}\Exp\left[\left\|\sqrt{C^{(n)}}^{-1} Y_j\right\|^3\right] \leq 
\frac{1}{n^{3/2}} \left\|\sqrt{C^{(n)}}^{-1}\right\|^3 b,
\] 
where $\|. \|$ denotes the operator norm. Here we have assumed that $C^{(n)}$ is invertible, which is certainly true for sufficiently large
$n$ since $C$ is invertible.
Since $\left\|\sqrt{C^{(n)}}^{-1}\right\|$ converges to the constant
$\gamma:=\|\sqrt{C}^{-1}\|$, we can 
bound $\sum_{j=1}^n \beta_j$ for all $n\geq n_0$ for sufficiently large  
$n_0$ by $(\gamma+\epsilon) b/\sqrt{n}$ with some fixed $\epsilon$.  
Let $Z_n$ be a Gaussian with covariance matrix $C^{(n)}$ and $Z$ be a Gaussian with 
covariance matrix $C$. 
Since the right hand side of
\eqref{eq:boundbentkus} converges to zero, 
we have
\begin{equation}\label{eq:convA}
\sup_{B \in {\cal C}} 
|P\{ S^{(n)}\in B \} - P\{Z_n \in B\}| \rightarrow 0. 
\end{equation}
Since $P\{ Z_n \in B \}$ converges to
$P\{ Z  \in B\}$ uniformly in $B$ (this is because the mapping of the covariance matrix to its Gaussian density is continuous at $C$ for the uniform norm topology on the mapping's codomain) \eqref{eq:convA} remains
true when $Z_n$ is replaced with $Z$.
Hence, $S^{(n)}$ converges in distribution  
to $Z$. 
 $\Box$
 
\vspace{0.3cm}

\noindent 
We now combine Lemma~\ref{lem:derived} and Lemma~\ref{lem:covasym} and obtain:
\begin{Lemma}[independence and Gaussianity of frequencies]\label{lem:indepGauss}
Given $k$ frequencies $\nu_1,\nu_2,\dots,\nu_k \in [-1/2,1/2]\setminus \{0\} $ with $|\nu_j|\neq |\nu_{j'}|$ for $j\neq j'$ and assume $\Exp[|A_j|^3]$ to be finite.  Then the sequence of random vectors
\[
(\hat{A}^n_{\nu_1},\hat{A}^n_{\nu_2},\dots,\hat{A}^n_{\nu_k}) \in \C^k
\]
converges in distribution to $(W_1,\dots,W_k)$ where $W_l$ are i.i.d.~random variables 
with distribution $G$ as in Theorem~\ref{thm:main}. 
\end{Lemma}

The proof is immediate after representing real and imaginary parts
of each $\hat{A}^n_{\nu_l}$ by an  $\R^2$-valued random variable as in Lemma~\ref{lem:covasym} and
 applying Lemma~\ref{lem:derived} to the random vector in $\R^{2k}$.
A uniform bound for $\Exp[\|S_j\|^3]$ follows easily from finiteness of
$\Exp[|A_j|^3]$. 

The fact that different $\hat{A}^n_\nu$ are asymptotically independent and identically distributed for different
$\nu$ has a very intuitive consequence: 
computing the Fourier sum $\hat{a}^n(\nu)$ for different frequencies and
one fixed instance $a=a_1,a_2,\dots$  
resembles the distribution of $\hat{A}^n_\nu$ for fixed $\nu$.
This suggests that the distribution of $\hat{A}^n_\nu(\omega)$ with $\nu$ uniformly chosen from $[-1/2,1/2]$ and fixed $\omega$ yields asymptotically also a Gaussian. To formally phrase
this idea we first need the following result:
\begin{Lemma}[distribution over a path]\label{lem:paths}
Let $(\tau,\Sigma_\tau,P_\tau)$ denote a
probability space and for each $t\in \tau$,
 let $(X^n_t)$ be a sequence of random vectors
 on the probability space $(\Omega,\Sigma,P_\Omega)$. 
Further, assume that
the map 
\[
X^n_.(\omega): t\mapsto X^n_t(\omega)
\]
is $\Sigma_\tau$-measurable for all $n\in \N$ and $\omega \in \Omega$ and thus defines
a random variable on $(\tau,\Sigma_\tau,P_\tau)$ 
whose distribution we denote by
$P_{X^n_.(\omega)}$. 

For every $k\in \N$ and $P_\tau^k$-almost all $k$-tuples $(t_1,\dots,t_k)$ let
the sequence of random vectors $(X^n_{t_1},\dots,X^n_{t_k})$ converge in distribution to $(Z_1,\dots,Z_k)$ for $n\to \infty$ where $Z_1,\dots,Z_k$ are i.i.d.~random variables on 
$(\Omega',\Sigma',P_{\Omega'})$
with distribution $P_Z$. 

Then the distance 
$d(X^n_.(\omega),P_Z)$ converges in probability to zero for any well-behaved 
pseudometric  $d$. More precisely,
the random variable
\[
\omega \mapsto d(X^n_.(\omega),P_Z)
\]
on $(\Omega,\Sigma,P_\Omega)$ 
converges to zero in probability.
\end{Lemma}

\noindent
Proof:  
We have to show that for every $\epsilon,\delta >0$ there is an $n_0$ such that
\[
P_\Omega \left\{ \omega\,\big|\, d(P_{X^n_.(\omega)}, P_Z) \geq \epsilon \right\} \leq \delta,
\]
for all $n\geq n_0$. 
For any
$k\in \N$, let $t_1,\dots,t_k$ be i.i.d.~drawn from $P_\tau$.
Since $d$ is well-behaved, the distance between 
$\hat{P}_{X^n_{t_1}(\omega),\dots,X^n_{t_k}(\omega)}$
and $P_{X^n_.(\omega)}$ converges  to zero in probability uniformly in $n$.
Thus, we can choose 
$k$ such that for all $n$
\begin{equation}\label{eq:firsterm}
d(P_{X^n_.(\omega)},\hat{P}_{X^n_{t_1}(\omega),\dots,X^n_{t_k}(\omega)}) \leq \epsilon/2
\end{equation}
holds with probability at least $1-\delta/3$, and that, at the same time,
\begin{equation}\label{eq:Zemp}
d(\hat{P}_{Z_1,\dots,Z_k},P_Z) \leq \epsilon/2
\end{equation} 
also holds with probability at least $1-\delta/3$. 
Using the triangle inequality for $d$ we obtain
\begin{equation}\label{eq:2terms}
d(P_{X^n_.(\omega)},P_Z) \leq d(P_{X^n_.(\omega)}, \hat{P}_{X^n_{t_1}(\omega),\dots,X^n_{t_k}(\omega)}) +
d(\hat{P}_{X^n_{t_1}(\omega),\dots,X^n_{t_k}(\omega)}, P_Z).
\end{equation}
The sequence of random vectors $(X_{t_1}^n,\dots,X_{t_k}^n)$ converge in
distribution to $(Z_1,\dots,Z_k)$ for each fixed $k$-tuple $(t_1,\dots,t_k)$. In other words, for any measurable set $B$ in $\R^k$  we have
\begin{equation}\label{eq:Bconv}
\lim_{n\to \infty} P_\Omega\{\omega \,|\, (X_1^n(\omega),\dots,X^n_k(\omega)) \in  B \}  = P_{\Omega'} \{ \omega'\,|\,(Z_1(\omega'),\dots,Z_k(\omega')) \in B\} .  
\end{equation}
Setting
\[
B:= \{z_1,\dots,z_k\,| d(\hat{P}_{z_1,\dots,z_k},P_Z) \leq \epsilon/2 \}, 
\]
we conclude from \eqref{eq:Bconv} that we can find an $n_0$ such that
\begin{equation}\label{eq:distreplace}
P_\Omega \{ \omega\,| d(\hat{P}_{X^n_{t_1}(\omega),\dots,X^n_{t_k}(\omega)},P_Z) \leq \epsilon/2 \}  \geq
P_{\Omega'} \{ \omega'\,|    d(\hat{P}_{Z_1(\omega'),\dots,Z_k(\omega')},P_Z) \leq \epsilon/2\} -\delta/3,
\end{equation}
for all $n\geq n_0$.
Using 
\eqref{eq:distreplace}
and
 \eqref{eq:Zemp} we can thus choose $n_0$ such that
for all $n\geq n_0$ 
\begin{equation}\label{eq:secondterm}
d(\hat{P}_{X^n_{t_1}(\omega),\dots,X^n_{t_k}(\omega)} , P_Z) \leq \epsilon/2
\end{equation}
with probability at most $1- 2 \delta/3$. 
Combining \eqref{eq:secondterm} with \eqref{eq:firsterm}
we have thus ensured that the right hand side of \eqref{eq:2terms}
is smaller than $\epsilon$ with probability at least $1-\delta$.  
$\Box$

\vspace{0.3cm}
\noindent
We are now able to prove 
Theorem~\ref{thm:main}
via Lemma~\ref{lem:paths}.
To this end,  let $(\tau,\Sigma_\tau,P_\tau)=
([-1/2/1/2],\cB,\lambda)$
 and set
\[
X^n_. (\omega) := ({\rm Re} \hat{A}^n_.(\omega), {\rm Im} \hat{A}^n_.(\omega)),
\]
 where ${\rm Re}$ and ${\rm Im}$ denote
real and imaginary part, respectively.    
If we then draw $d$ frequencies $\nu_1,\dots,\nu_d$ 
for arbitrarily large $d$, we satisfy $P^k_\tau$ -almost surely the condition $0 \neq |\nu_j| \neq |\nu_{j'}| \neq 0$
required by Lemma~\ref{lem:indepGauss}. 
Thus, the sequence of random vectors $(\hat{A}^n_{\nu_1},\dots,\hat{A}^n_{\nu_k})$ converges in distribution to 
$(W_1,\dots,W_k)$ where $W_j$ are 
distributed according to $G$.
Therefore, the random variable
\[
\omega \mapsto d(P_{\hat{A}^n_.(\omega)},G)
\]
converges to zero in probability due to Lemma~\ref{lem:paths}.


\end{document}